\documentclass[11pt,reqno]{amsart}
\usepackage{graphicx,color}
\usepackage{amsmath}
\usepackage{amssymb}
\usepackage{amstext}
\usepackage[T1]{fontenc}
\usepackage[latin1]{inputenc}
\usepackage{mathrsfs}
\date{}

\newtheorem{rem}{Remark}[section]
\newtheorem{thm}{Theorem}[section]

\newtheorem{lem}[thm]{Lemma}
\newtheorem{prop}[thm]{Proposition}
\newtheorem{defi}[thm]{Definition}

\newenvironment{prooft}{ \rm\textit{Proof of Theorem
\ref{thm4.1}.}\par}{${}$\hfill$\square$}

\numberwithin{equation}{section}
\begin{document}

\title[\textbf{Existence and relaxation theorems }]
{\bf Existence and relaxation theorems for differential inclusions involving time dependent maximal monotone operators }%
\author{Amira. Makhlouf}
\address{Amira Makhlouf\\ Laboratoire LAOTI, Universit\'e Mohammed Seddik Benyahia-Jijel,
Alg\'erie} \email{amira.makhlouf18@gmail.com,
a\_makhlouf@univ-jijel.dz}\author{Dalila. Azzam-Laouir}
\address{Dalila Azzam-Laouir\\ Laboratoire LAOTI, Universit\'e Mohammed Seddik Benyahia-Jijel,
Alg\'erie} \email{laouir.dalila@gmail.com, dalilalaouir@univ-jijel.dz}
\author{Charles. Castaing}
\address{Charles.Castaing\\ IMAG,
Universit\'e Montpellier, CNRS, Montpellier Cedex 5, France}
\email{charles.castaing@gmail.com}

\maketitle

\begin{abstract}
In this paper, we consider evolution problems involving time
dependent maximal monotone operators in  Hilbert spaces. Existence
and relaxation theorems are proved.
\\\vspace{0.4cm}\\
\textit{Keywords.} Absolutely continuous variation, convexified
problem, extreme points, fixed  point, maximal monotone operator,
pseudo-distance, perturbation, relaxation, weak norm.
\\\textit{Mathematics Subject Classification (2010).}
34A60; 34B15; 47H10.
\end{abstract}
\section{Introduction}
Let $H$ be a separable Hilbert space and $I=[0,T]$ $(T>0)$ be an
interval of $\mathbb{R}$. In this paper, we establish    existence
of solutions of problems$$(\mathcal{P}_{f,F})\quad
\left\{\begin{array}{l}
      -\dot u(t)\in
A(t)u(t)+f(t,u(t))+F(t,u(t)), \ a.e.\,  t\in I, \\
u(t)\in D(A(t)),\ a.e.\,  t\in I, \\
                                                                       u(0)=u_0\in D(A(0)),    \\
                                                                     \end{array}
\right.$$
 and $$(\mathcal{P}_{f,ext (F)})\quad  \left\{\begin{array}{l}
      -\dot u(t)\in
A(t)u(t)+f(t,u(t))+ext (F(t,u(t))), \ a.e.\,   t\in I \\u(t)\in
D(A(t)),\  a.e.\,  t\in I ,\\
                                                                       u(0)=u_0\in D(A(0)),     \\
                                                                     \end{array}
\right.$$  and we  also consider the problem
$$(\mathcal{P}_{f,\overline{co}(F)})\quad  \left\{\begin{array}{lc}
      -\dot u(t)\in
A(t)u(t)+f(t,u(t))+\overline{co}(F(t,u(t))), \ a.e. \,   t\in I
\\u(t)\in D(A(t)),\ a.e.\,  t\in I, \\
                                                                       u(0)=u_0\in D(A(0)),    \\
                                                                     \end{array}
\right.$$ where, $f:I\times H\to H$ is separately integrable on $I$
and separately Lipschitz on $H$, $F:I\times H\rightrightarrows H$ is
a set-valued map with compact values, and, for all $t\in I$, $A(t)$
is a maximal monotone operator of $H$ and $D(A(t))$ its domain. The
dependence $t\mapsto A(t)$ is absolutely continuous   with respect
to the pseudo-distance, $dis(\cdot, \cdot)$, introduced by
Vladimirov \cite{V} (see relation \eqref{2.3}).

Moreover, by taking  $H=\mathbb{R}^d$, the $d$-dimensional Euclidean
space, we prove under a suitable hypothesis on $F$, that the
solutions set of problem $(\mathcal{P}_{f,ext ( F)})$ is dense in
the solutions set of the problem $(\mathcal{P}_{f,F})$ and the
solutions set of problem $(\mathcal{P}_{f,F})$ is dense in the
solutions set of the problem $(\mathcal{P}_{f,\overline{co}(F)})$
(relaxation theorems).

Problem $(\mathcal{P}_{f,F})$ has been studied in many papers with
different type of perturbations, and where the variation of the time
dependent maximal monotone operator $A(t)$ is absolutely continuous,
Lipschitz or BVC (continuous with bounded variation), we refer to
\cite{ABCM,ABCM2,AB,ACM,KM,Tols1,Tols2,V}. These works constitute
new development of existing ones in the literature dealing with
differential inclusions governed by fixed maximal monotone operators
(not depending in the time), see for instance \cite{AD,Ba,Ba2,B} and
their references.

Existence and relaxation  problems have been considered by many
authors and exist in the literature
\cite{AP,AMT,AM,C,GT,DP,DP2,F,T,W} among other.

 Existence of solutions for differential inclusions
with extreme points of a set-valued map as right-hand side, has been
studied by many authors, see \cite{Br,C,Ch1,Ch2,DP,DP2,TT2,TT}. In
some of these works we find application to density theorems. The
main key in the proof of such problems is the existence of
continuous selections with values in the extreme points of a
set-valued map. Some extensions to second order differential
inclusions were established by  Avgerinos et al. in \cite{AP}, Gomaa
in \cite{Go}, Ibrahim et al. in \cite{IG}, and recently  by Azzam et
al. in \cite{AM} for a second order problem governed by the
subdifferential of a convex function.

To the best of our knowledge, this work is the first investigating
relaxation problems with maximal monotone operator, which is an
extension of the results by Tolstonogov \cite{Tols1} and Tolstonogov
et al. \cite{Tols2}  dealing with relaxation  problems with
subdifferential of convex functions and sweeping processes, since
sweeping process is a differential inclusion governed by $N_{C(t)}$;
the normal cone of a moving closed convex set $C(t)$, which is a
particular case of maximal monotone operators, and since if we take
$A(t)=N_{C(t)}$, we have
$$ dis(A(t), A(s))=\mathcal{H}(C(t), C(s))$$
where $\mathcal {H}$ denotes the Hausdorff distance between closed
sets.

The plan of the paper is the following. In the next section, we give
notation and  some preliminary results necessary for our
investigations. In section 3, we  present two results on the
existence of solutions of problems $(\mathcal{P}_{f,F})$ and
$(\mathcal{P}_{f,ext(F)})$ by following the ideas of the proof of
the results in \cite{AB} and \cite{AM}. Finally, we prove relaxation
theorems in section 4, by using  techniques in  \cite{AP},
\cite{AM}, \cite{T} and \cite{Tols1}.
\section{Preliminaries}

Let $I=[0,T]$ $(T>0)$ be an interval of $\mathbb{R}$, and let $X$ be
a separable Banach space with norm $\|\cdot\|_X$ and $X'$ its
topological dual. We introduce the following notation: $Id_X$ the
identity mapping of $X$, $\mathcal{L}(I )$ the Lebesgue
$\sigma$-algebra on $I$, $\mu= dt$  the Lebesgue measure on $I$,
$\mathcal{B}(X)$  the $\sigma$-algebra of Borel subsets of $X$,
$\overline{B}(x_0, r)$ the closed ball of $X$ of center $x_0$ and
radius $r > 0$ and $\overline{B}_X$ the closed unit ball. \\For
closed subsets $S_1 $ and $S_2 $ of $X$, the Hausdorff distance
between them is defined
by$$\mathcal{H}(S_1,S_2)=\max\left\{\sup_{x\in S_2}d(x,S_1),\
\sup_{x\in S_1}d(x,S_2)\right\},$$ where $d(x,S)=\inf
\left\{\|x-y\|_X:\ y\in S\right\}$, for $S\subset X$.
\\ One defines the (possibly empty)
set of nearest points of $x$ in $S$ by $$Proj_S(x) := \big\{u \in S
:\ d(x,S) = \|x- u\|_X\big\}.$$When $Proj_S(x)$ contains one and
only one point $u_0$, we will write $u_0 = proj_S(x)$.

 For any subset $S$ of $X$, $S^c$ stands for the complement of $S$,  $co(S)$  for its  convex hull
 and
$\overline{co}(S)$ for its closed convex hull. If $X$ is a
finite-dimensional space and $S$ is compact, then so is $co (S)$.

 We denote by $\mathrm{C}(I,X)$ the Banach space of all continuous
mappings $u : I\to X$, endowed with the sup norm
$\|u\|_C=\max\limits_{t\in I}\|u(t)\|_X$.\\ By $\mathrm{L}^p(I,X)$
($1\le p \le +\infty$), we denote the Banach space of all
equivalence classes of  measurable mappings $u:I\to X$, equipped
with its standard norm $\|\cdot \|_p$.\\
By $\mathrm{L}^p_\sigma(I,X)$ ($1\le p < +\infty$), we denote the
space $\mathrm{L}^p(I,X)$ equipped with the weak norm
$\|.\|_{\sigma}$ (see \cite{G}) which is defined by
\begin{equation}\label{2.2} \|u\| _{\sigma}= \max_{ t\in I}
\left\|\int^t_0 u(s)ds\right\|_X,\ \hbox{for } \ u\in
\mathrm{L}^p(I,X).
\end{equation}For $p\in [1,+\infty[$, $\mathrm{W}^{1,p}(I,X)$ is the space of the
absolutely continuous mappings $u:I\to X$ such that $\dot u\in
\mathrm{L}^p(I,X)$.

We denote by $\sigma$-$X$ the space $X$ endowed with the weak
topology $\sigma(X,X')$. The same notation is used for subsets of
$X$. Otherwise, we assume that the space $X$ and its subsets are
endowed with the strong (normed) topology.

Let $X$,  $Y$ be two topological spaces. A set-valued map $F
:X\rightrightarrows Y$ with nonempty  values   is said to be upper
semicontinuous  (resp. lower semicontinuous) if  the inverse image
$F^{-1}(C) = \{x \in X : F (x) \cap C \not =\emptyset\}$ is closed
(resp. open)  for any closed (resp. open) subset $C$ of $Y$.

Let $(\Omega,\Sigma)$  be a measurable space and $X$ be a separable
Banach space. Let  $F :\Omega\rightrightarrows X$ be a set-valued
map with nonempty values. $F$ is said to be measurable or
$\Sigma$-measurable if $F^{-1}(C) \in \Sigma$ for each closed subset
$C$ of $ X$. When we say that $F$ is graph measurable, we mean that
$gph(F)=\big\{(\omega,x)\in \Omega\times X:\ x\in F(\omega)\big\}\in
\Sigma\otimes \mathcal{B}(X)$. If $F$ is a measurable set-valued map
with nonempty and closed values then $F$ is graph measurable, and we
have equivalence if $(\Omega,\Sigma,\mu)$ is a $\sigma$-finite and
complete measure space. Also, if  $F$ is a measurable set-valued map
with nonempty and closed values then,  the function $\omega\mapsto
d(x, F(\omega))$, $x\in X$, and the set-valued maps $\omega\mapsto
co(F(\omega))$ and $\omega\mapsto \overline{co}(F(\omega))$ are
measurable. See \cite{CV}.

{\lem\cite{CV}\label{lem1} Let $(\Omega,\Sigma)$ be a measurable
space and $X$ be a separable Banach space. Let $F : \Omega \times
X\rightrightarrows X$ be a measurable set-valued map and  $u :
\Omega\to X$ be a measurable mapping. Then, the set-valued map
$\omega\mapsto F(\omega, u(\omega))$ is measurable.}

{\prop \label{prop1}\cite{BC} Let $X$ be a metric space and $Y$ be a
Banach space. Let $F: X\rightrightarrows \mathrm{L}^1(I,Y)$ be a
lower semicontinuous set-valued map with closed and decomposable
values.
\\Assume that $g:X\to \mathrm{L}^1(I,Y)$ and $\phi: X\to \mathrm{L}^1(I,\mathbb{R})$ are continuous functions such that, for every $x\in X$, the set
 $$\Psi(x)=\Big\{ u\in F(x):\ \|u(t)-g(x)(t)\|_Y< \varphi(x)(t) \ \text{a.e. on} \ I\Big\}$$ is nonempty.
 Then the set-valued map $\Psi:X \rightrightarrows \mathrm{L}^1(I,Y)$
 is lower semicontinuous
  and has  decomposable values.}

Here the notion "decomposable" means that if $S$ is a measurable
subset of $I$, and $g_1$, $ g_2 \in  F(x)$, then $g_1
1\!\mathrm{I}_S+g_2 1\!\mathrm{I}_{S^c}  \in F(x)$ (see \cite{Ba}),
where $1\!\mathrm{I}_S$ denotes the characteristic function of $S$.

 Let $H$ be a
real separable  Hilbert space with the  scalar product
$\langle\cdot,\cdot\rangle$ and the associated norm $\|\cdot\|$.
Given a set-valued map $F:I\rightrightarrows H$ with nonempty
values and
 $p\in[1,+\infty]$, we denote  by $S^p_F$  the set of selections of
$F$ which belong to  $L^p(I,H)$, i.e., $S^p_F =\{f(\cdot) \in L^p(I,
H) : f (t) \in F(t) \;\hbox{a.e. on}\; I\}$. In general this set may
be empty. However,  for a graph measurable set-valued map $F$,  the
set $S^p_F$ is nonempty if and only if $t\mapsto \inf\{\|v\| : v \in
F(t)\} \in L^p(I,\mathbb{R}_+)$ (see  \cite[Lemma 2.3.2]{HP1}).

    Now, we give the definition and some properties of maximal
monotone operators. We refer the reader to \cite{Ba,B,V2}  for this
concept and details. A set-valued operator of $H$ is a set-valued
map from $H$ to $H$, its domain, range and graph will be denoted
here, and hereafter by $D(A)$, $R(A)$ and $gph(A)$, these sets are
defined by
$$D(A)=\{x\in H:\ Ax\neq \emptyset\},\quad R(A) =\{y\in H:\ \exists x\in D(A),\ y\in
Ax\},$$$$gph(A)=\{(x,y)\in D(A)\times H:\ y\in Ax\}.$$ We say that
 $A:D(A)\subset H \to 2^H$ is monotone if
$\langle y_1 -y_2, x_1- x_2\rangle\ge0$ whenever $(x_i,y_i)\in
gph(A),\ i = 1, 2$.  We say that a monotone operator $A$ is maximal
if $gph(A)$ is not contained properly in any other monotone
operator, that is by Minty's Theorem, for all $\lambda> 0$, $R(Id_H
+\lambda A) = H$. If $A$ is a maximal monotone operator then, for
every $x\in D(A)$, $Ax$ is a nonempty, closed and convex set. We
denote the projection of $0$ into $Ax$, $proj_{Ax}(0)$, by $A^0(x)$.

 Now, for $\lambda > 0$ we
define the resolvent of $A$,  $J_{\lambda}= (Id_H + \lambda
A)^{-1}$, and the Yosida approximation of $A$, $A_{\lambda} =
\frac{1}{ \lambda} (Id_H- J_{\lambda}).$ These operators are defined
on all of $H$.

Let $A : D(A)\subset H \to 2^H$ and $B : D(B)\subset H\to2^H$ be two
maximal monotone operators, then we denote by $dis(A,\ B)$ the
pseudo-distance between $A$ and $B$ defined by (see \cite{V})
\begin{equation}\label{2.3}dis(A,\ B) = \sup
\left\{\frac{\langle y - \hat{y}, \hat{x} - x\rangle}{ 1 + \|y\| +
\|\hat{y}\| }:\ (x,y) \in gph(A),(\hat x, \hat y)\in gph(B)
\right\}.
\end{equation}
\begin{lem}{\rm\cite{KM}}
Let $A$ be a maximal monotone operator. If $x\in \overline{D(A)}$
and $ y \in H$ are such that $$\langle A^0(\eta) -y, \eta -
x\rangle\ge  0,\quad \hbox{for}\ \eta\in D(A),$$ then $x\in D(A)$
and $y\in Ax$.
\end{lem}

 In what
follows, we recall a series of results, which we will use in our
proofs.

{\prop\label{prop2.1.1} \cite{HP1}  Suppose that $X$ is a separable
metric space, and $Y$ is a Banach space. Suppose that $F:
X\rightrightarrows \mathrm{L}^p(I,Y) \ (p\ge1)$ is a lower
semicontinuous set-valued map with nonempty, closed and decomposable
values. Then $F$ has a continuous selection, i.e., there exists a
continuous mapping $f: X \to \mathrm{L}^p(I,Y)$, satisfying $f(x)\in
F(x)$ for all $x\in X$.}

{\prop \label{prop2.1} Let $F : I\times H\rightrightarrows H$ be a
set-valued map with nonempty and closed values. Assume
that\begin{enumerate}
\item $gph(F)\in \mathcal{L}(I)\otimes \mathcal{B}(H)\otimes
\mathcal{B}(H)$;
    \item  the set-valued
map $x\mapsto F(t, x)$  is lower semicontinuous  a.e. on $I$;
    \item there exists a nonnegative function $m\in
\mathrm{L}^2(I,\mathbb{R})$ such that
\begin{equation}\label{2.5} \sup\{\|v\|:\ v\in F(t, u(t))\}\le
m(t),\ \ a.e.\, \textmd{on}\; I\quad\forall  u(\cdot)\in {C}(I,H).
\end{equation}\end{enumerate} Then, there exists a continuous mapping $g :{C}(I,H)\to
\mathrm{L}^2(I,H)$ such that $$g(u)(t)\in F(t, u(t)), \ \ a.e.\,
\textmd{on}\; I\quad\forall u(\cdot)\in {C}(I,H) .$$}
\begin{proof}
Define $N:{C}(I,H)\rightrightarrows \mathrm{L}^2(I,H)$ by
$$N(u)=\big\{h\in \mathrm{L}^2(I,H):\ h(t)\in F(t,u(t))\ \ a.e.\, \textmd{on}\; I\big\}\quad\forall u\in {C}(I,H).$$  Using the same techniques
 in the proof of Proposition III.2.7 in \cite{HP}, we can show
that $N$ is lower semicontinuous with nonempty, closed and
decomposable values. Hence, by Proposition \ref{prop2.1.1},  there
exists a continuous mapping $g :{C}(I,H)\to \mathrm{L}^2(I,H)$ such
that
$$g(u)\in N(u), \quad\hbox{for all}\quad u(\cdot)\in {C}(I,H),$$
that is,
$$g(u)(t)\in F(t, u(t)), \ \ a.e.\, \textmd{on}\; I\quad\forall u(\cdot)\in
{C}(I,H).  $$
\end{proof}

{\prop\cite[Proposition 6.10]{TT}\label{prop2.2} Assume that $p\in
[1,+\infty[$. Let $X$ be a separable Banach space,
$\mathcal{K}\subset C(I,X)$ be a compact set, and let $F : I\times
X\rightrightarrows X$ be a set-valued map with nonempty and compact
values. Assume that\begin{enumerate}
\item  the set-valued
map $t\mapsto F(t, x)$  is measurable   for all $x\in X$;
    \item  the set-valued
map $x\mapsto  F(t, x) $  is Hausdorff continuous  almost everywhere
on $I$;
    \item there exists a nonnegative function $m\in
\mathrm{L}^p(I,\mathbb{R})$ such that \eqref{2.5} is satisfied.
\end{enumerate}      Then, for any continuous mapping  $g :\mathcal{K}\to
\mathrm{L}^p(I,X)$ such that $$g(u)(t)\in \overline{co} (F(t,
u(t))),\ \hbox{a.e. on}\ I \quad\hbox{for all}\quad u(\cdot)\in
\mathcal{K}   ,$$ and for any $\varepsilon
> 0$, there exists a continuous mapping $g_\varepsilon :\mathcal{K}\to \mathrm{L}^p(I,X)$ such that
$$g_\varepsilon(u)(t)\in F(t, u(t))\quad a.e.\,
\textmd{on}\; I,$$ and $$\|g(u) - g_\varepsilon(u)\|_{\sigma}<
\varepsilon, \quad\forall u(\cdot)\in \mathcal{K}.$$}{\defi\cite{AM}
Let $K$ be a subset of a linear space $X $. We say that a point
$x\in X$ is an extreme point of $K$, and we write $x\in ext(K)$, if
for all $x_1, x_2\in K$ and $\lambda \in]0, 1[,\ x = \lambda x_1 +
(1 -\lambda)x_2$ implies that $x_1 = x_2$.} {\thm\cite{AM} Let $X$
be a topological vector space whose dual space $X'$ separates points
of $X$. Then every nonempty compact convex subset $K$ of $X$ has at
least one extreme point, i.e., $ext(K) \neq \emptyset$.} {\prop
\cite[Proposition 8.1.]{TT2}\label{prop2.1.21}  Assume that $p\in
[1,+\infty[$. Let $X$ be a separable Banach space,
$\mathcal{K}\subset C(I,X)$ be a compact set, and let $F : I\times
X\rightrightarrows X$ be a set-valued map with nonempty, convex and
weakly compact values. Assume that\begin{enumerate}
\item  the set-valued
map $t\mapsto F(t, x)$  is measurable   for all $x\in X$;
    \item  the set-valued
map $x\mapsto F(t, x)$  is   continuous  almost everywhere on $I$;
    \item there exists a nonnegative function $m\in
\mathrm{L}^p(I,\mathbb{R})$ such that \eqref{2.5} is satisfied.
\end{enumerate}      Then, there exists a  continuous mapping $g :\mathcal{K}\to
\mathrm{L}^p(I,X)$ such that $$g(u)(t)\in ext( F(t, u(t))),\  \
a.e.\, \textmd{on}\; I \quad\forall u(\cdot)\in \mathcal{K}.$$}
{\prop\cite[Proposition 6.9]{TT}\label{prop2.1.2} Assume that $p\in
[1,+\infty[$. Let $X$ be a separable Banach space,
$\mathcal{K}\subset C(I,X)$ be a compact set, and let $F : I\times
X\rightrightarrows X$ be a set-valued map with nonempty, convex and
weakly compact values. Assume that\begin{enumerate}
\item  the set-valued
map $t\mapsto F(t, x)$  is measurable   for all $x\in X$;
    \item  the set-valued
map $x\mapsto F(t, x)$  is   continuous  almost everywhere on $I$;
    \item there exists a nonnegative function $m\in
\mathrm{L}^p(I,\mathbb{R})$ such that \eqref{2.5} is satisfied.
\end{enumerate}      Then, for any continuous mapping $g :\mathcal{K}\to
\mathrm{L}^p(I,X)$ such that $$g(u)(t)\in  F(t, u(t)),\  \ a.e.\,
\textmd{on}\; I \quad\forall u(\cdot)\in \mathcal{K},$$  and for any
$\varepsilon
> 0$, there exists a continuous mapping $g_\varepsilon :\mathcal{K}\to \mathrm{L}^p(I,X)$ such that
$$g_\varepsilon(u)(t)\in ext( F(t, u(t))) \ \ a.e.\,
\textmd{on}\; I,$$ and $$\|g(u) - g_\varepsilon(u)\|_{\sigma}<
\varepsilon, \quad\forall u(\cdot)\in \mathcal{K}.$$} {\lem{\rm
\cite{T,G}}\label{lem2.3} Let $X = \mathbb{R}^d$, $m \in
\mathrm{L}^1(I,\mathbb{R})$ be  a nonnegative function, and $$G =
\{f \in \mathrm{L}^1(I,\mathbb{R}^d):\ \|f(t)\|\le m(t) \ \ a.e.\,
\textmd{on}\; I\}. $$ Then the topologies of the spaces
$\sigma-\mathrm{L}^1(I,\mathbb{R}^d)$ and
$\mathrm{L}^1_\sigma(I,\mathbb{R}^d)$ coincide on  $G$ and
consequently, the set $G$ is convex, metrizable and compact in
$\mathrm{L}^1_\sigma (I,\mathbb{R}^d)$. } \vskip2mm We close this
section by the  following Gronwall's Lemma. {\lem{\rm
\cite{B}}\label{lem2.4} Let $m \in \mathrm{L}^1(I,\mathbb{R})$ such
that $m (t)\ge 0$ a.e. on $]0, T[$ and let $\alpha$ be a nonnegative
constant. Let $\phi : I\to \mathbb{R}$ be a continuous function
satisfying $$\phi(t) \le \alpha +\int^t_0 m(s)\phi(s)ds, \;\;
\forall\ t \in I,$$ then
$$|\phi(t)|\le \alpha e^{\int^t_0 m(s)ds},\;\; \forall\ t\in I.$$

    \section{Existence results}
   \rm

For the statement of our theorems of  this section we have to assume
the following hypotheses.
  $(\mathcal{H}_A^1)$ There exists a function
$\beta\in
    W^{1,2}(I,\mathbb{R})$ which is nonnegative on $[0,T[$ and
    nondecreasing  such that
    $$dis(A(t),A(s))\le |\beta(t)-\beta(s)|,\;\; \forall\, t,s\in
    I.$$
    $(\mathcal{H}_A^2)$ There exists $c\ge 0$ such that
    $$\|A^0(t,x)\|\le c(1+\|x\|),\;\;\forall\, t\in I \; \text{and}\  x\in
    D(A(t)).$$
    $(\mathcal{H}_A^3)$ For each $x\in H$, the mapping $t\mapsto
    A_1(t)x$ is measurable.\\
    $(\mathcal{H}_A^4)$  For every $t\in I$,  $D(A(t))$  is relatively ball
    compact.\\
    $(\mathcal{H}_f^1)$ For every $R>0$, there is a nonnegative
    integrable function $\lambda_R$ such that, for all $t\in I$
    $$\|f(t,x)-f(t,y)\|\le  \lambda_R(t)\|x-y\|,\;\;\forall x,y\in
    \overline{B}(0,R).$$
    $(\mathcal{H}_f^2)$ There exists a nonnegative real number $M$ such that
    $$\|f(t,x)\|\le M (1+\|x\|), \; \;\forall\ t\in I,\ x\in H.$$
    $(\mathcal{H}_F^1)$ $F$ is $(\mathcal{L}(I)\otimes
    \mathcal{B}(H))$-measurable.\\
     $(\mathcal{H}_F^2)$ For all $t\in I$, the set-valued map $x\mapsto F(t,x) $ is lower semicontinuous on
     $H$.\\
      $(\mathcal{H}_F^3)$ There exists a nonnegative function $m\in \mathrm{L}^2(I,\mathbb{R})$ such that
      $$F(t,x)\subset m(t)\overline{B}_{H},\ \ \text{for all}\ (t,x)\in I\times
    H.$$
$(\mathcal{H}_F^4)$ The set-valued map $t\mapsto F(t, x)$  is
measurable   for all $x\in H$.\\
    $(\mathcal{H}_F^5)$ The set-valued
map $x\mapsto F(t, x)$  is  Hausdorff continuous  a.e. on $I$.
\vskip2mm
    For the proof of our theorems we will need the following results
    from \cite{ABCM} and \cite{AB}. {\thm \label{thm3.1} Let for every $t\in I$,
    $A(t):D(A(t))\subset H \to 2^H$ be a maximal monotone operator
    satisfying $(\mathcal{H}_A^1)$ and $(\mathcal{H}_A^2)$. Let
    $f:I\times H\to H$ be such that for every $x\in H$
    $f(\cdot,x)$ is measurable on $I$. Suppose also that $f$
    satisfies $(\mathcal{H}_f^1)$ and $(\mathcal{H}_f^2)$. Let
    $h:I\to H$ be a  mapping in  $\mathrm{L^2}(I,H)$. Then, for all
    $u_0\in D(A(0))$, the problem $$(\mathcal{P}_{f,h})\quad  \left\{\begin{array}{l}
      -\dot u(t)\in                                                    A(t)u(t)+f(t,u(t))+h(t), \ a.e. \ t\in I, \\
      u(t)\in
D(A(t)),\  a.e.\  t\in I ,\\
                                                                       u(0)=u_0\in D(A(0)), \\
                                                                     \end{array}
\right.$$ admits a unique  absolutely continuous solution $u$.
Moreover,$$\|\dot u(t)\|\le K(1+\dot
\beta(t))+(1+K)\|h(t)\|\quad\hbox{ a.e.}\ t\in I$$ for some
nonnegative real constant $K=K(u_0,c,M,h,\beta,T)$.} }
{\prop\label{prop3.1} Let for every $t\in I$, $A(t):D(A(t))\subset H
\to 2^H$ be a maximal monotone operator
    satisfying $(\mathcal{H}_A^1)$, $(\mathcal{H}_A^2)$ and
    $(\mathcal{H}_A^3)$.\\Let $(u_n)$ and $(v_n)$ be sequences in  $\mathrm{L^2}(I,H)$ satisfying \begin{enumerate}
        \item $v_n(t)\in A(t) u_n(t) $ for all $n\in \mathbb{N}$ and
        almost every $t\in I$;
        \item  $(u_n)$ converges strongly to $u\in
        \mathrm{L^2}(I,H)$;
        \item  $(v_n)$ converges weakly to $v\in \mathrm{L^2}(I,H)$.
    \end{enumerate}         Then, we have $u(t)\in D(A(t))$ and  $v(t)\in A(t)u(t)$ for almost every $t\in I$.}
\vskip2mm
 Now, we present our first existence theorem, it can be
seen as a particular case of Theorem 3.5. in \cite{AB}, since in
\cite{AB} the perturbation $F$ is mixed semicontinuous, that is, for
every $t\in I$, at each $x\in H$ such that $F(t,x)$ is convex the
set-valued map $F(t,\cdot)$ is upper semicontinuous and where
$F(t,x)$ is not convex, the set-valued map $F(t,\cdot)$ is lower
semicontinuous on some neighborhood of $x$. {\thm\label{thm3.2} Let
for every $t\in I$,
    $A(t):D(A(t))\subset H \to 2^H$ be a maximal monotone operator
    satisfying $(\mathcal{H}_A^1)$, $(\mathcal{H}_A^2)$, $(\mathcal{H}_A^3)$ and $(\mathcal{H}_A^4)$.  Let
    $f:I\times H\to H$ be such that for every $x\in H$
    $f(\cdot,x)$ is measurable on $I$. Suppose also that $f$
    satisfies $(\mathcal{H}_f^1)$ and $(\mathcal{H}_f^2)$. Let
    $F:I\times H\rightrightarrows H$ be a set-valued map with
    nonempty and  compact values satisfying
    $(\mathcal{H}_F^1)$, $(\mathcal{H}_F^2)$ and $(\mathcal{H}_F^3)$. Then, for all
    $u_0\in D(A(0))$,  problem $$(\mathcal{P}_{f,F})\quad  \left\{\begin{array}{l}
      -\dot u(t)\in                                                    A(t)u(t)+f(t,u(t))+F(t,u(t)),\ a.e.  \ t\in I ,\\
      u(t)\in
D(A(t)),\  a.e.\  t\in I ,\\
                                                                       u(0)=u_0\in D(A(0)),  \\
                                                                     \end{array}
\right.$$ admits an absolutely continuous (a.c) solution  $u$.
Moreover,
$$\|\dot u(t)\|\le K(1+\dot \beta(t))+(1+K)m(t)=:\gamma(t)\quad \hbox{ a.e.
}\ t\in I$$ where $K=K(\|u_0\|, T, c, M, m,\beta)$.}

\begin{proof}
Using Proposition \ref{prop2.1}, we obtain a continuous mapping
$g:{C}(I,H) \to \mathrm{L}^2 (I,H)$ satisfying
\begin{equation}\label{3.1.1}
g(u)(t)\in F(t, u(t)), \; \;a.e.\, \textmd{on}\; I \quad \forall
u(\cdot)\in {C}(I,H)  .
\end{equation} Let us consider the set $$W=\big\{h\in
\mathrm{L}^2(I,H):\ \|h(t)\|\le m(t) \;a.e.\, t\in I\big\}.$$ It is
clear that $W$ is convex and by Banach-Alaoglu-Bourbaki's theorem,
it is a weakly compact subset of $\mathrm{L}^2(I,H)$. Let us
consider the set $$\Lambda=\big\{u_h:\ u_h\ \text{is the unique a. c
solution of }\ (\mathcal{P}_{f,h}),\ h\in W\big\}.$$ According to
Theorem \ref{thm3.1},  $\Lambda$ is nonempty and for each $h\in W$,
we have
\begin{equation}\label{a} \|\dot u_h(t) \|\le\gamma(t) \ \
\hbox{a.e. on } \ I ,
\end{equation} where $\gamma(t)=K(1+\dot
\beta(t))+(1+K)m(t)$ for all $t\in I$ and $K$ is the constant in
Theorem \ref{thm3.1} with $m$ instead of $\|h(\cdot)\|$. Then, for
each $h\in W$ and for all $s,t\in I$ with $s\le t$ we have
$$\|u_h(t)-u_h(s)\|=\left\|\int_s^t \dot u_h(\tau) d\tau \right\|\le
\int_s^t \gamma(\tau)d\tau.$$ As $\dot \beta,m\in \mathrm{L}^2
(I,\mathbb{R})$, we get that $\gamma\in\mathrm{L}^2 (I,\mathbb{R})$
and so $\gamma \in \mathrm{L}^1 (I,\mathbb{R})$. We conclude that
$\Lambda$ is equicontinuous in ${C}(I,H)$.\\On the other hand, we
have for all $t\in I$,  $\Lambda(t)=:\{u_h(t):\ u_h\in \Lambda\}
\subset D(A(t))$ and for all $h\in W$
$$\|u_h(t)\| \le \|u_0\|+\int_0^t\|\dot u_h(s)\|ds\le
\|u_0\|+\int_0^t\gamma(s)ds\le \|u_0\|+\|\gamma\|_1=: R.$$ Hence,
$\Lambda(t)\subset D(A(t))\cap\ \overline{B}(0,R)$. Then, it is
relatively compact according to $(\mathcal{H}_A^4)$. By virtue of
Arzelà-Ascoli's theorem, $\Lambda$ is relatively compact in
${C}(I,H)$.

Let us set $\widehat{\Lambda}=\overline{co}(\Lambda)$. It is clear
that $\widehat{\Lambda}$ is convex and compact in ${C}(I,H)$.

In the following, consider the mapping $\phi:W\to \widehat{\Lambda}$
defined by $\phi(h)=u_h$ where $u_h$ is the unique absolutely
continuous solution of $(\mathcal{P}_{f,h})$ and let us prove that
$\phi$ is continuous from $W$ endowed with the weak topology into
$\widehat{\Lambda}$. For this purpose, let $(h_n)$ be a sequence of
$W$ converging to $h$ with respect to the weak topology. Then,
$(\phi(h_n))=(u_{h_n})$ is a sequence of $\widehat{\Lambda}$ such
that  for every $n\in \mathbb{N}$
$$-\dot u_{h_n}(t)\in A(t)u_{h_n}(t)+f(t,u_{h_n}(t))+h_n(t)\
\text{a.e. on }\ I,$$ and $u_{h_n}(0)=u_0$. Since
$\widehat{\Lambda}$ is compact, we can extract from $(u_{h_n})$ a
subsequence, that we do not relabel, which converges to some mapping
$w\in\widehat{\Lambda}$. Also, since $\|\dot u_{h_n}(t)\|\le
\gamma(t) $ a.e. $t\in I$, for all $n\in \mathbb{N}$, and since
$\gamma\in \mathrm{L}^2(I,\mathbb{R})$, we conclude that $(\dot
u_{h_n})$ is bounded in $\mathrm{L}^2(I,H)$, so we can extract a
subsequence, not relabeled, which converges weakly in
$\mathrm{L}^2(I,H)$ to some mapping  $\zeta\in \mathrm{L}^2(I,H)$.
For all $t\in I$ and $x\in H$,  $ 1\!\mathrm{I}_{[0,t]}x\in
\mathrm{L}^2(I,H)$. Then, we have for each fixed $t$
 \begin{multline*}
   \Big\langle \lim_{n\to \infty}\int_0^t\dot u_{h_n}(s)ds,x\Big\rangle =\lim_{n\to\infty}\Big\langle\int_0^t \dot
 u_{h_n}(s)ds,x\Big\rangle = \lim_{n\to\infty}\int_0^T\langle \dot
 u_{h_n}(s), 1\!\mathrm{I}_{[0,t]}(s)x\rangle ds
\\=\lim_{n\to\infty}\langle \dot
 u_{h_n}, 1\!\mathrm{I}_{[0,t]}x\rangle
= \langle \zeta, 1\!\mathrm{I}_{[0,t]}x\rangle
   =  \int_0^T\langle \zeta(s), 1\!\mathrm{I}_{[0,t]}(s)x\rangle ds
 =   \Big\langle\int_0^t \zeta(s) ds,x\Big\rangle .
 \end{multline*}So that, $\lim\limits_{n\to\infty}\int_0^t \dot
 u_{h_n}(s)ds= \int_0^t \zeta(s)$. Hence, $$w(t)=\lim_{n\to\infty}u_{h_n}(t)=u_0+\lim_{n\to\infty}\int_0^t\dot
 u_{h_n}(s)ds=u_0+\int_0^t\zeta(s)ds.$$This implies that $\dot
 w=\zeta$ for almost every $t\in I$. We conclude that  $(\dot
 u_{h_n})$ converges weakly to $\dot w$ in
 $\mathrm{L}^2(I,H)$.

 On the other hand, as $(u_{h_n}(t))\subset\overline{B}(0,R)$ for all $t\in
 I$, by
 $(\mathcal{H}^2_f)$ and the uniform convergence of $(u_{h_n})$ to
 $w$, we get, for each $t\in I$, $$\lim_{n\to \infty}
 f(t,u_{h_n}(t))=f(t,w(t)),$$ and by $(\mathcal{H}^1_f)$
 $$\|f(t,u_{h_n}(t))\|\le M(1+\|u_{h_n}(t)\|)\le M(1+R).$$ By the
 dominated convergence theorem, $(f(\cdot,u_{h_n}(\cdot)))$ converges in
 $\mathrm{L}^2(I,H)$ to $f(\cdot,w(\cdot))$ and then, it converges
 weakly in  $\mathrm{L}^2(I,H)$ to this  limit. \\Since,  for every $n\in \mathbb{N}$ $$-\dot
u_{h_n}(t)-f(t,u_{h_n}(t))-h_n(t)\in A(t)u_{h_n}(t)\ \text{a.e. on
}\ I,$$ and  $\Big(\dot
u_{h_n}(\cdot)-f(\cdot,u_{h_n}(\cdot))-h_n(\cdot)\Big)$ converges
weakly in $\mathrm{L}^2(I,H)$ to $\dot
w(\cdot)-f(\cdot,w(\cdot))-h(\cdot)$, we conclude by Proposition
\ref{prop3.1}, that $w(t)\in D(A(t))$ and $$-\dot w(t)\in
A(t)w(t)+f(t,w(t))+h(t)\ \text{a.e.}\ t\in I,$$ with
$w(0)=\lim\limits_{n\to \infty}u_{h_n}(0)=u_0$, that is $w$ is the
unique absolutely continuous solution to   problem
$(\mathcal{P}_{f,h})$ and consequently, $w=u_h\in\widehat{\Lambda}$,
so that $(\phi(h_n))$ converges to $\phi(h)$ in $\widehat{\Lambda}$.
This shows the required  continuity of $\phi$.

Now, note that, by \eqref{3.1.1},   $g(\widehat{\Lambda})\subset W
$, so let $\widehat g$ be  the restriction of $g$ to
$\widehat{\Lambda}$ and let us define  $\psi=\phi\circ \widehat
g:\widehat{\Lambda}\to \widehat{\Lambda}$ and  prove that it is also
 a continuous mapping.\\Indeed, let $(u_n)\subset\widehat{\Lambda}$ and assume that
$(u_n)$ converges to the map $\bar{u}\in\widehat{\Lambda}$. Then
$(\widehat{g}(u_n))$ converges to $\widehat{g}(\bar{u})$ in
$\mathrm{L}^2(I,H)$, and so $(\widehat{g}(u_n))$ converges weakly to
$\widehat{g}(\bar{u}) $ in $\mathrm{L}^2(I,H)$. By what preceds
$(\phi(\widehat{g}(u_n)))$ converges to
$\phi(\widehat{g}(\bar{u}))$. Then $(\psi(u_n))$ converges
to $\psi(\bar{u})$ in $\widehat{\Lambda}$. Whence $\psi$ is continuous.\\
An application of Shauder's fixed point theorem asserts the
existence of  some element $u\in \widehat{\Lambda}$ such that
$u=\psi(u)=(\phi\circ \widehat{g})(u)=\phi(g(u))$, that is
$$-\dot u(t)\in A(t)u(t)+f(t,u(t))+g(u)(t)\ \text{a.e. on}\ I,$$ and
$u(0)=u_0\in D(A(0))$. Since $g(u)(t)\in F(t,u(t))$ a.e. on $I$, we
conclude that
$$\left\{\begin{array}{lcl}
      -\dot u(t)\in                                                    A(t)u(t)+f(t,u(t))+ F(t,u(t)),\;  a.e.\,   t\in I\\
      u(t)\in D(A(t)), \;  a.e.\,   t\in I\\
                                                                       u(0)=u_0\in D(A(0)).
                                                                     \end{array}
\right.$$ So that $u\in\mathrm{W}^{1,2}(I,H)$ is a solution of
$(\mathcal{P}_{f,F})$. Furthermore,  $$\|\dot u(t)\|\le
\gamma(t)\;\; a.e. t\in I.$$
\end{proof}
Next, we give our second existence theorem. {\thm\label{thm3.3} Let
for every $t\in I$,
    $A(t):D(A(t))\subset H \to 2^H$ be a maximal monotone operator
    satisfying $(\mathcal{H}_A^1)$, $(\mathcal{H}_A^2)$, $(\mathcal{H}_A^3)$ and $(\mathcal{H}_A^4)$.  Let
    $f:I\times H\to H$ be such that for every $x\in H$
    $f(\cdot,x)$ is measurable on $I$. Suppose also that $f$
    satisfies $(\mathcal{H}_f^1)$ and $(\mathcal{H}_f^2)$. Let
    $F:I\times H\rightrightarrows H$ be a set-valued map with
    nonempty, convex and compact values satisfying
    $(\mathcal{H}_F^3)$, $(\mathcal{H}_F^4)$ and $(\mathcal{H}_F^5)$. Then, for all
    $u_0\in D(A(0))$, the problem $$(\mathcal{P}_{f,ext (F)})\quad  \left\{\begin{array}{lcl}
      -\dot u(t)\in                                                    A(t)u(t)+f(t,u(t))+ext (F(t,u(t))), \; a.e. \, t\in I \\
       u(t)\in D(A(t)), \;  a.e.\,   t\in I\\
                                                                       u(0)=u_0 \\
                                                                     \end{array}
\right.$$ admits an absolutely continuous solution $u$. Moreover,
$$\|\dot u(t)\|\le K(1+\dot \beta(t))+(1+K)m(t)=:\gamma(t)\;\; a.e.
\, t\in I,$$ where $K=K(\|u_0\|, T, c, M, m,\beta)$.}

\begin{proof}
We use the same notations of the proof of Theorem \ref{thm3.1}, that
is, we  consider the convex weakly compact subset of
$\mathrm{L}^2(I,H)$
$$W=\big\{h\in \mathrm{L}^2(I,H):\ \|h(t)\|\le m(t) \;\; a.e.\, t\in
I\big\},$$ and the relatively compact subset of ${C}(I,H)$
$$\Lambda=\big\{u_h:\ u_h\ \text{is the unique a.c solution of }\
(\mathcal{P}_{f,h}),\ h\in W\big\}$$   and we  set
$\widehat{\Lambda}=\overline{co}(\Lambda)$, which is a convex
compact subset of  ${C}(I,H)$.

By Proposition \ref{prop2.1.21}, there exists a continuous mapping
$\widehat{g}:\widehat{\Lambda}\to \mathrm{L}^2(I,H)$ such that
$$\widehat{g}(u)(t)\in ext (F(t,u(t))),\ \ a.e.\,\textmd{ on}   \; I\quad \forall u\in \widehat{\Lambda},$$ that is, for all $u\in
\widehat{\Lambda}$, $\widehat{g}(u)\in S^2_{ext(F(\cdot,u(\cdot))}$
and then, $\widehat{g}(\widehat{\Lambda})\subset W$.

As in the proof of Theorem \ref{thm3.2}, we define the mapping
$\psi=\phi\circ \widehat{g}: \widehat{\Lambda}\to \widehat{\Lambda}$
and we apply Schauder's fixed point theorem to find a point
$u\in\widehat{\Lambda}$ such that $u=(\phi\circ \widehat{g})(u)$,
that is
$$-\dot u(t)\in A(t)u(t)+f(t,u(t))+\widehat{g}(u)(t)\ \ \text{a.e. on}\ I,$$
with $u(t)\in D(A(t))$ a.e. and $u(0)=u_0\in D(A(0))$, and since
$\widehat{g}(u)(t)\in ext (F(t,u(t)))$ we get
$$-\dot
u(t)\in A(t)u(t)+f(t,u(t))+ext (F(t,u(t)))\ \ \text{a.e. on}\ I$$
with $u(t)\in D(A(t))$ a.e. and $u(0)=u_0\in D(A(0))$, which means
that $u$ is an absolutely continuous solution of
$(\mathcal{P}_{f,ext( F)}). $ Furthermore, as $\phi(W)\subset
\Lambda$, we get $u\in\Lambda$ and so $\|\dot u(t)\|\le \gamma(t)$
a.e. on $I$.
\end{proof}

{\rem   The hypotheses of Theorem \ref{thm3.2} ensure the existence
of solutions of  problem $(\mathcal{P}_{f,\overline{co}(F)})$.}

{\rem If $H=\mathbb{R}^d$,  hypothesis $(\mathcal{H}_A^4)$ can be
omitted in both theorems \ref{thm3.2} and \ref{thm3.3}.}
\section{Relaxation theorems} In this section we will establish some relaxation theorems related to our existence results in the previous section.
For this purpose we have to take $H=\mathbb{R}^d$ and to assume the
following hypothesis on $F$.\\
$(\mathcal{H}_F^6)$ There exists a nonnegative function $k\in
\mathrm{L}^1(I,\mathbb{R})$ such that
$$\mathcal{H}(F(t,x),
    F(t,y))\le k(t)\|x-y\|,  \quad\text{for all} \ (t,x,y)\in I\times H\times
    H.$$
    Let us consider the convexified problem
$$(\mathcal{P}_{f,\overline{co}F})\quad \left\{\begin{array}{l}
      -\dot u(t)\in                                                    A(t)u(t)+f(t,u(t))+\overline{co}(F(t,u(t))), \ a.e.   \ t\in I, \\
      u(t)\in D(A(t)),\ a.e. \ t\in I,\\
                                                                       u(0)=u_0\in D(A(0)). \\
                                                                     \end{array}
\right.$$

{\thm\label{thm4.1}   Let for every $t\in I$,
    $A(t):D(A(t))\subset H\to 2^H$ be a maximal monotone operator
    satisfying $(\mathcal{H}_A^1)$, $(\mathcal{H}_A^2)$ and  $(\mathcal{H}_A^3)$.  Let
    $f:I\times H\to H$ be such that for every $x\in H$
    $f(\cdot,x)$ is measurable on $I$. Suppose also that $f$
    satisfies $(\mathcal{H}_f^1)$ and $(\mathcal{H}_f^2)$. Let
    $F:I\times H\rightrightarrows H$ be a set-valued map
    with
    nonempty and  compact values satisfying $(\mathcal{H}_F^1)$, $(\mathcal{H}_F^3)$ and $(\mathcal{H}_{F}^6)$. Then the solution
     set $S(\mathcal{P}_{f,F})$ of the problem $(\mathcal{P}_{f,F})$ is dense in the solution set  $S(\mathcal{P}_{f,\overline{co}(F)})$ of the problem
     $(\mathcal{P}_{f,\overline{co}(F)})$ with respect
to the topology of uniform convergence.} {\rem Since, in Theorem
\ref{thm4.1}, $F$ has compact values and $H=\mathbb{R}^d$, the
set-valued map $(t,x)\mapsto co(F(t,x))$ has also compact values, so
that $\overline{co}(F(t,x))=co(F(t,x))$ for every $(t,x)\in I\times
H$.}

      \noindent \begin{prooft} First, note that the sets $W$,
      $\Lambda$ and $\widehat{\Lambda}$ are the same in the proof of
      Theorem \ref{thm3.2}.
Let $u(\cdot)\in S(\mathcal{P}_{f, {co}(F)})$, then, there exists
$z\in S^2_{{co}(F(\cdot,u(\cdot)))}$ such that
$$(\mathcal{P}_{f,z})\quad \left\{\begin{array}{l}
      -\dot u(t)-f(t,u(t))-z(t)\in                                                    A(t)u(t),\ a.e.  \ t\in I, \\
      u(t)\in D(A(t)),\ a.e.\ t\in I,\\
                                                                       u(0)=u_0\in D(A(0)).  \\
                                                                     \end{array}
\right.$$ Notice   that    $\|\dot u(t)\|\le \gamma(t)$, a.e. $ t\in
I$, so that, \begin{equation}\label{4.1.3} \|u(t)\|\le
\|u_0\|+\|\gamma\|_1=: R,\ \ \forall  t\in I.
\end{equation}
Let
$\varepsilon>0$ and $w\in C(I,H)$, and let us define the set-valued
map $\Phi_\varepsilon: I\rightrightarrows H$ by
$$\Phi_\varepsilon(t)=\Big\{y\in {co}(F(t,w(t))):\
\|z(t)-y\|<\varepsilon+d(z(t),{co}( F(t,w(t))))\Big\}.$$ Evidently,
$\Phi_\varepsilon(t)\neq \emptyset$ for all $t\in I$, and we have
\begin{eqnarray*}
  gph(\Phi_\varepsilon) &=&\Big\{ (t,y)\in I\times H:\ y\in \Phi_\varepsilon (t)\Big\}\\
  &=& \Big\{ (t,y)\in I\times H:\
\|z(t)-y\|<\varepsilon+d(z(t), {co}( F(t,w(t))))\Big\}\cap gph\big(
{co}(F(\cdot,w(\cdot)))\big).
\end{eqnarray*} By hypothesis $(\mathcal{H}_F^1)$ and Lemma \ref{lem1}, the mapping
$t\mapsto d(z(t), {co}(F(t,w(t))))$ and the set-valued map $t\mapsto
 {co} (F(t,w(t)))$ are $\mathcal{L}(I)$-measurable. Then $
gph(\Phi_\varepsilon)\in \mathcal{L}(I)\otimes \mathcal{B}(H)$.\\
Apply the measurable selection theorem (see \cite[Theorem
III.6.]{CV}) to obtain a measurable map $v:I\to H$ such that
$v(t)\in \Phi_{\varepsilon}(t)$, for all $t\in I$.

 Next, we define the set valued map
 $\Psi_\varepsilon:\widehat{\Lambda}\rightrightarrows\mathrm{L}^1(I,H)$
 by $$\Psi_\varepsilon(w)=\Big\{g\in
 S^1_{{{co}(F(\cdot,w(\cdot)))}}:\
 \|z(t)-g(t)\|<\varepsilon+d(z(t),{co}(F(t,w(t))))\
 \text{a.e. on}\ I\Big\},$$ where $\widehat{\Lambda}$ is the convex
 compact set defined in the proof of Theorem \ref{thm3.2}.\\
By what precedes, we have for all $w\in \widehat{\Lambda}$, $v\in
\Psi_\varepsilon(w)$, so that $\Psi_\varepsilon(w) \neq \emptyset$.
Hence, from Proposition \ref{prop1}
 $w\mapsto \Psi_\varepsilon(w)$ is lower semicontinuous
with decomposable values, so is $w\mapsto
\overline{\Psi_\varepsilon(w)}$. By Proposition \ref{prop2.1.1}, we
get a continuous mapping $g_\varepsilon:
\widehat{\Lambda}\to\mathrm{L}^1(I,H)$ such that
$g_\varepsilon(w)\in \overline{\Psi_\varepsilon(w)}$ for all $w\in
\widehat{\Lambda}$. Then, for all $w\in \widehat{\Lambda}$,
$g_\varepsilon(w)(t)\in {co}(F(t,w(t)))$ a.e. on $I$ and
$$\|z(t)-g_\varepsilon(w)(t)\|\le\varepsilon+d(z(t),{co}(F(t,w(t))))\
 \text{a.e. on}\ I.$$\\ An application of
Proposition \ref{prop2.2}, gives us a continuous mapping
$\varphi_\varepsilon:\widehat{\Lambda}\to \mathrm{L}^1(I,H)$ such
that
$$\varphi_\varepsilon(w)(t)\in F(t,w(t))\ \ a.e.\, \textmd{on}\; I,$$ and
$$\|g_\varepsilon(w)-\varphi_\varepsilon(w)\|_{\sigma}<\varepsilon\ \ \forall  w\in \widehat{\Lambda}.$$
We take in the following a sequence
$(\varepsilon_n)$ of nonnegative real numbers which decreases to $0$
as $n\to \infty$. Then, by the arguments above, for each $n\in
\mathbb{N}$, we have mappings $g_{\varepsilon_n}$ and
$\varphi_{\varepsilon_n}$ satisfying for all $w\in\hat{\Lambda}$,
$g_{\varepsilon_n}(w)\in \overline{\Psi_\varepsilon(w)}$,
$\varphi_{\varepsilon_n}\in S^1_{F(\cdot,w(\cdot))}$ and
$\|g_{\varepsilon_n}(w)-\varphi_{\varepsilon_n}(w)\|_{\sigma}<\varepsilon_n$.
Furthermore, $\varphi_{\varepsilon_n}(\widehat{\Lambda})\subset W$
and $g_{\varepsilon_n}(\widehat{\Lambda})\subset W$.

To be more clear, we will index our sequences and sets by $n$
instead of $\varepsilon_n$.

Let us consider, for every $n$, the set-valued map
$\Gamma_n:W\rightrightarrows \mathrm{L}^1(I,H)$ defined by
$$\Gamma_n(h)=\{\varphi_n(\phi(h))\}\ \ \forall h\in W. $$
So, $\Gamma_n$ has nonempty, convex  and closed values and
$\Gamma_n(W)\subset W$, that is $\Gamma_n$ maps $W$ into itself. Let
us prove that $\Gamma_n$ is upper semicontinuous from $W$ into
itself endowed with the weak topology of $\mathrm{L}^2(I,H)$. As $W$
is weakly compact in $\mathrm{L}^2(I,H)$, it suffices to show that
$gph(\Gamma_n)$ is sequentially weakly closed in $W\times W$. For
this end, let $(h_k,l_k) $ be a sequence in $gph(\Gamma_n)$ such
that $(h_k,l_k) $ converges weakly to $(h,l)\in W\times W$, then for
each $k\in\mathbb{N}$, $l_k\in\Gamma_n(h_k)$, that is,
$l_k=\varphi_n(\phi(h_k))$. As $\phi$ is continuous from $W$ endowed
with the weak topology into $\widehat{\Lambda}$, $(\phi(h_k))$
converges to $\phi(h)$ in $\widehat{\Lambda}$, then, by the
continuity of $\varphi_n$, $(\varphi_n(\phi(h_k)))$ converges to
$\varphi_n(\phi(h))$ in $\mathrm{L}^1(I,H)$, so
$(\varphi_n(\phi(h_k)))$ converges weakly to $\varphi_n(\phi(h))$ in
$\mathrm{L}^1(I,H)$. But $(l_k)$ converges weakly to $l$ in
$\mathrm{L}^2(I,H)$, and so in $\mathrm{L}^1(I,H)$. Hence,
$l=\varphi_n(\phi(h))\in \Gamma_n(h)$, that is $(h,l)\in
gph(\Gamma_n)$. This shows that $gph(\Gamma_n)$ is sequentially
weakly closed in $W\times W$ and hence we get the upper
semicontinuity of $\Gamma_n$.

An application of Kakutani-Ky Fan fixed point theorem gives some
element $h_n\in W$ such that $h_n\in \Gamma_n(h_n)$. We put
$u^*_n=\phi(h_n)$, then   $u^*_n\in \widehat{\Lambda}$
   and
   \begin{equation}\label{4.1.1}
-\dot u^*_n(t)\in A(t)u^*_n(t)+f(t,u^*_n(t))+\varphi_n(u^*_n)(t)\
\text{a.e. on}\ I,
\end{equation}
    with $u^*_n(0)=u_0$ and $\|\dot  u^*_n(t)\|\le \gamma(t)$ a.e.
   on $I$, so that
   \begin{equation}\label{4.1.2}
\|u^*_n(t)\|\le \|u_0\|+\|\gamma\|_1=R \ \ \forall t\in
   I.
\end{equation}
    Since $(u^*_n)\subset \widehat{\Lambda}$
and $\widehat{\Lambda}$ is compact, we can extract a subsequence,
that we do not relabel, which converges to some mapping $u^*\in
\widehat{\Lambda}$.

 Whence, referring to $(\mathcal{P}_{f,z})$,
\eqref{4.1.3}, \eqref{4.1.1} and \eqref{4.1.2} and using
 the monotonicity of $A(t)$ as well as $(\mathcal{H}_f^1)$, we have, for almost all $t\in I$ \begin{eqnarray*}
   &&\frac 1 2 \frac d{dt} \|u^*_n(t)-u(t)\|^2 = \langle\dot u^*_n(t)-\dot u(t),u^*_n(t)-u(t)\rangle  \\
    &=&\big\langle\dot u^*_n(t)+f(t,u^*_n(t))+\varphi_n(u^*_n)(t)-\dot u(t)-f(t,u(t))-z(t),u^*_n(t)-u(t)
    \big\rangle\\& +&\big\langle z(t)-\varphi_n(u^*_n)(t),u^*_n(t)-u(t)\big\rangle+\big\langle
     f(t,u(t))-f(t,u^*_n(t)) ,u^*_n(t)-u(t)\big\rangle  \\
    &\le& \big\langle z(t)-\varphi_n(u^*_n)(t),u^*_n(t)-u(t)\big\rangle +\big\langle
     f(t,u(t))-f(t,u^*_n(t)) ,u^*_n(t)-u(t)\big\rangle   \\
    &\le&\big\langle z(t)-g_n(u^*_n)(t),u^*_n(t)-u(t)\big\rangle +\big\langle
    g_n(u^*_n)(t)-\varphi_n(u^*_n)(t),u^*_n(t)-u(t)\big\rangle\\&+&\lambda_R(t)\|u^*_n(t)-u(t)\|^2.
 \end{eqnarray*}
 Integrating between $0$ and $t$ and taking into account that
 $u^*_n(0)=u(0)=u_0$, we obtain
 \begin{multline}\label{4.1}
   \frac 1 2\|u^*_n(t)-u(t)\|^2 \le \int_0^t\langle z(s)-g_n(u^*_n)(s),u^*_n(s)-u(s)\rangle ds\\
    +\int_0^t\langle g_n(u^*_n)(s)-\varphi_n(u^*_n)(s),u^*_n(s)-u(s)\rangle
    ds+\int_0^t\lambda_R(s)\|u^*_n(s)-u(s)\|^2ds.
 \end{multline} We know that $\|g_n(u^*_n)-\varphi_n(u^*_n)\|_{\sigma}<\varepsilon_n$,
  then as $n\to \infty$, $(g_n(u^*_n)-\varphi_n(u^*_n))$ converges to $0$ in
 $\mathrm{L}^1_\sigma(I,H)$. Since $g_n(u^*_n)(t)\in {co}(F(t,u^*_n(t)))$ a.e.\ on $I$
 and $\varphi_n(u^*_n)(t)\in F(t,u^*_n(t))$  a.e.\, on $I$, it is clear
 that, for each $n$,
  $g_n(u^*_n),\varphi_n(u^*_n)\in D=\big\{h\in \mathrm{L}^1(I,H):\ \|h(t)\|\le m(t) \ \ a.e.\;
  on\,
 I\big \}$.
Then, by Lemma \ref{lem2.3}, $(g_n(u^*_n)-\varphi_n(u^*_n))$
converges weakly  on $\mathrm{L}^1(I,H)$ to $0$, i.e., for all
$\xi\in \mathrm{L}^{\infty}(I,H)$ $\lim_{n\to \infty} \langle
g_n(u^*_n)-\varphi_n(u^*_n),\xi\rangle =0$, in particular, for $t\in
I$ and $\xi=1\!\mathrm{I}_{[0,t]}(u^*-u)$  we have
\begin{equation}\label{4.2}
\lim_{n\to \infty} \langle
g_n(u^*_n)-\varphi_n(u^*_n),1\!\mathrm{I}_{[0,t]}(u^*-u)\rangle =0.
\end{equation}
So, for all $t\in I$,
 $$\lim_{n\to \infty} \int _0^t\langle g_n(u^*_n)(s)-\varphi_n(u^*_n)(s),u^*(s)-u(s)\rangle ds =0.$$
 On the other hand, since $(g_n(u^*_n))$ and $(\varphi_n(u^*_n))$ are
 bounded  in $\mathrm{L}^2(I,H)$, and so in $\mathrm{L}^1(I,H)$, we have
 \begin{equation}\label{4.3}
\int _0^t\langle
g_n(u^*_n)(s)-\varphi_n(u^*_n)(s),u^*_n(s)-u^*(s)\rangle ds \le
2\|m\|_1\|u^*_n-u^*\|_C\underset{n\to \infty}\to 0.
\end{equation}
By \eqref{4.2} and \eqref{4.3} we get
$$ \lim_{n\to \infty} \int _0^t\langle g_n(u^*_n)(s)-\varphi_n(u^*_n)(s),u^*_n(s)-u(s)\rangle ds =0.$$
Furthermore, we have by $(\mathcal{H}^6_F)$ and the fact that $z\in
S^2_{{co}(F(\cdot,u(\cdot)))}= S^1_{{co}(F(\cdot,u(\cdot)))} $ and
$g_n(u^*_n)\in \overline{\Psi_{n}(u^*_n)}$,
\begin{eqnarray*}
 && \int_0^t\langle z(s)-g_n(u^*_n)(s),u^*_n(s)-u(s)\rangle ds\\  &\le& \int_0^t\| z(s)-g_n(u^*_n)(s)\|\|u^*_n(s)-u(s)\|ds\\
   &\le& \int_0^t\bigg(\varepsilon_n+d\big(z(s), {co}(F(s,u^*_n(s)))\big)\bigg) \|u^*_n(s)-u(s)\|ds\\
   &\le&  \int_0^t\bigg(\varepsilon_n+\mathcal{H}\big( {co}(F(s,u(s))), {co}(F(s,u^*_n(s)))\big)\bigg) \|u^*_n(s)-u(s)\|ds\\
    &\le&  \int_0^t\bigg(\varepsilon_n+\mathcal{H}\big(F(s,u(s)),F(s,u^*_n(s))\big)\bigg) \|u^*_n(s)-u(s)\|ds\\
   &\le&\int_0^t\bigg(\varepsilon_n+k(s)\|u(s)-u^*_n(s)\|\bigg) \|u^*_n(s)-u(s)\|ds.
\end{eqnarray*}
Consequently, letting $n\to \infty$ in \eqref{4.1}, we get
$$\frac 1 2\|u^*(t)-u(t)\|^2\le \int_0^t(k(s)+\lambda_R(s)) \|u^*(s)-u(s)\|^2ds.$$ Then, by Lemma \ref{lem2.4}, we obtain that $$\|u^*(t)-u(t)\|\le 0, \
\text{for all}\ t\in I,$$ that is $u^* =u$. Consequently, $u^*_n\to
u$ in $\mathrm{C}(I, H)$ as $n\to \infty$ with $(u^*_n)\subset
S(\mathcal{P}_{f,F})$
 and this proves that $S(\mathcal{P}_{f,{co}(F)})\subset \overline{S(\mathcal{P}_{f,F})}$.

 To finish the proof, we need to show that  $S(\mathcal{P}_{f,{co}F})$ is closed. To this end, let $(u_n)\subset
 S(\mathcal{P}_{f,{co}(F)})$ and assume that $(u_n)$ converges uniformely to $u\in \mathrm{C}(I, H)$.
  Then $$-\dot u_n(t)\in A(t)u_n(t)+f(t,u_n(t))+z_n(t)\ \ \text{a.e. on } \ I,$$
with $u_n(0)=u_0$ and $z_n\in S^2_{ {co}(F(\cdot,u_n(\cdot)))}$,
i.e., $z_n(t)\in {co}(F(t,u_n(t)))$ a.e. on $I$ for all $n$.

From hypothesis $(\mathcal{H}_F^3)$ we have  $(z_n)\subset W$, and
as $W$  is weakly compact in $\mathrm{L}^2(I,H)$,  we can extract a
subsequence converging weakly to some mapping $z\in W$. Then, by
using the same arguments as in the proof oft he continuity of
$\phi$, we have that
$$-\dot u(t)\in A(t)u(t)+f(t,u(t))+z(t)\ \ \text{a.e. on } \ I,$$
with $u(0)=u_0$. Moreover,
 since $ {co}(F(t,\cdot))$ is Hausdorff continuous with nonempty, convex and compact values in $H$, we
 get by Theorem VI-4 in \cite{CV}
 $$z(t)\in {co}(F(t,u(t)))\quad \text{a.e. on} \ I,$$
 that is,$$-\dot u(t)\in A(t)u(t)+ f(t,u(t))+ co(F(t,u(t)))\quad \text{a.e. on} \ I.$$
 Consequently, $u\in S(\mathcal{P}_{f,{co}(F)})$. Hence $S(\mathcal{P}_{f,{co}(F)})$  is closed in $\mathrm{C}(I, H)$.\\Since
 ${S(\mathcal{P}_{f,F})}\subset S(\mathcal{P}_{f,{co}(F)})$, we conclude that $S(\mathcal{P}_{f,{co}(F)})= \overline{S(\mathcal{P}_{f,F})}$,
 and this completes our proof. \end{prooft}
\vskip2mm
 Next, we give our second relaxation theorem.
 {\thm   Let for every $t\in I$,
    $A(t):D(A(t))\subset H \to 2^{H}$ be a maximal monotone operator
    satisfying $(\mathcal{H}_A^1)$, $(\mathcal{H}_A^2)$ and  $(\mathcal{H}_A^3)$.  Let
    $f:I\times H\to H$ be such that for every $x\in H$
    $f(\cdot,x)$ is measurable on $I$. Suppose also that $f$
    satisfies $(\mathcal{H}_f^1)$ and $(\mathcal{H}_f^2)$. Let
    $F:I\times H\rightrightarrows H$ be a set-valued map
    with
    nonempty, convex and compact values satisfying $(\mathcal{H}_F^1)$, $(\mathcal{H}_F^3)$ and $(\mathcal{H}_{F}^6)$. Then the solution
     set $S(\mathcal{P}_{f,ext (F)})$ of the problem $(\mathcal{P}_{f,ext (F)})$ is dense in the solution set  $S(\mathcal{P}_{f,F})$
     of the problem $(\mathcal{P}_{f,F})$
      with respect to the topology of uniform convergence.}

      \begin{proof} We consider the sets $W$,
      $\Lambda$ and $\widehat{\Lambda}$ are the same in the proof of
      Theorem \ref{thm3.2}.

Let $u(\cdot)\in S(\mathcal{P}_{f,F})$, then, there exists $z\in
S^2_{F(\cdot,u(\cdot))}$ such that
$$(\mathcal{P}_{f,z})\quad \left\{\begin{array}{lc}
      -\dot u(t)-f(t,u(t))-z(t)\in                                                    A(t)u(t), & \text{a.e.}\  t\in I \\
                                                                       u(0)=u_0\in D(A(0)). &   \\
                                                                     \end{array}
\right.$$Notice  that    $\|\dot u(t)\|\le \gamma(t)$, a.e. $ t\in
I$, so that, \begin{equation}\label{4.1.32} \|u(t)\|\le
\|u_0\|+\|\gamma\|_1= R,\ \ \text{ for all} \ t\in I.
\end{equation}\\Let
$\varepsilon>0$ and $w\in C(I,H)$, and let us define the
set-valued map $\Phi_\varepsilon: I\rightrightarrows H$ by
$$\Phi_\varepsilon(t)=\Big\{v\in F(t,w(t)):\
\|z(t)-v\|<\varepsilon+d(z(t), F(t,w(t)))\;\text{a.e. on}\
I\Big\}.$$ Then, as in the proof of Theorem \ref{thm4.1}, there
exists  a measurable map $v:I\to H$ such that $v(t)\in
\Phi_{\varepsilon}(t)$, for all $t\in I$.

 Hence, we define the set valued map
 $\Psi_\varepsilon:\widehat{\Lambda}\rightrightarrows\mathrm{L}^1(I,H)$
 by $$\Psi_\varepsilon(w)=\Big\{g\in
 S^1_{F(\cdot,w(\cdot))}:\
 \|z(t)-g(t)\|<\varepsilon+d(z(t),F(t,w(t)))\
 \text{a.e. on}\ I\Big\}.$$
We have for all $w\in \widehat{\Lambda}$, $\Psi_\varepsilon(w)\neq
\emptyset$ and so, arguing as in the proof of Theorem \ref{thm4.1},
we get a continuous mapping $g_\varepsilon:
\widehat{\Lambda}\to\mathrm{L}^1(I,H)$ such that
$g_\varepsilon(w)\in \overline{\Psi_\varepsilon(w)}$ for all $w\in
\widehat{\Lambda}$, that is,  $g_\varepsilon(w)(t)\in F(t,w(t))$
a.e. on $I$ for all $w\in \widehat{\Lambda}$ and
$\|z(t)-g_\varepsilon(w)(t)\|\le\varepsilon+d(z(t),F(t,w(t)))$, for
all $w\in \widehat{\Lambda}$ and a.e. on $I$.\\ An application of
Proposition \ref{prop2.1.2}, gives us a continuous mapping
$\varphi_\varepsilon:\widehat{\Lambda}\to \mathrm{L}^1(I,H)$ such
that
$$\varphi_\varepsilon(w)(t)\in ext (F(t,w(t)))\ \ \text{a.e. on}\
I,$$ and
$$\|g_\varepsilon(w)-\varphi_\varepsilon(w)\|_{\sigma}<\varepsilon,\ \text{for
all} \ w\in \widehat{\Lambda}.$$We take in the following a sequence
$(\varepsilon_n)$ of nonnegative real numbers which decreases to $0$
as $n\to \infty$. Then, for each $n\in \mathbb{N}$, we have mappings
$g_{\varepsilon_n}$ and $\varphi_{\varepsilon_n}$ satisfying for all
$w\in\hat{\Lambda}$, $g_{\varepsilon_n}(w)\in
\overline{\Psi_\varepsilon(w)}$, $\varphi_{\varepsilon_n}(w)\in
S^1_{F(\cdot,w(\cdot))}$ and
$\|g_{\varepsilon_n}(w)-\varphi_{\varepsilon_n}(w)\|_{\sigma}<\varepsilon_n$.
Furthermore, $\varphi_{\varepsilon_n}(\widehat{\Lambda})\subset W$
and $g_{\varepsilon_n}(\widehat{\Lambda})\subset W$.

As above, we will index our sequences and sets by $n$
instead of $\varepsilon_n$.

As in the proof of  Theorem \ref{thm4.1}, for each $n\in\mathbb{N}$,
we define  the upper semicontinuous set-valued map
$\Gamma_n:W\rightrightarrows \mathrm{L}^1(I,H)$ by
 $$\Gamma_n(h)=\{\varphi_n(\phi(h))\},\ \ \text{for all}\ h\in W, $$ and we apply  Kakutani-Ky Fan fixed point theorem to
get some element  $h_n\in W$ such that $h_n\in \Gamma_n(h_n)$. We
put $u^*_n=\phi(h_n)$, then   $u^*_n\in \widehat{\Lambda}$
   and  \begin{equation}\label{4.1.12}
-\dot u^*_n(t)\in A(t)u^*_n(t)+f(t,u^*_n(t))+\varphi_n(u^*_n)(t)\
\text{a.e. on}\ I,
\end{equation}  with $u^*_n(0)=u_0$ and $\|\dot  u^*_n(t)\|\le \gamma(t)$ a.e.
   on $I$, so that
   \begin{equation}\label{4.1.22}
\|u^*_n(t)\|\le \|u_0\|+\|\gamma\|_1 =R\ \ \text {for all}\ t\in
   I.
\end{equation}

     Since $(u^*_n)\subset \widehat{\Lambda}$
and $\widehat{\Lambda}$ is compact, we can extract a subsequence,
that we do not relabel, which converges to some mapping $u^*\in
\widehat{\Lambda}$.

Whence, referring to $(\mathcal{P}_{f,z})$,
\eqref{4.1.32}, \eqref{4.1.12} and \eqref{4.1.22} and using
 the monotonicity of $A(t)$ as well as $(\mathcal{H}_f^1)$, we have, for almost all $t\in I$\begin{eqnarray*}
   &&\frac 1 2 \frac d{dt} \|u^*_n(t)-u(t)\|^2 = \langle\dot u^*_n(t)-\dot u(t),u^*_n(t)-u(t)\rangle  \\
    &=&\langle\dot u^*_n(t)+f(t,u^*_n(t))+\varphi_n(u^*_n)(t)-\dot u(t)-f(t,u(t))-z(t),u^*_n(t)-u(t)\rangle\\
    & +&\langle z(t)-\varphi_n(u^*_n)(t),u^*_n(t)-u(t)\rangle+\langle
     f(t,u(t))-f(t,u^*_n(t)) ,u^*_n(t)-u(t)\rangle  \\
    &\le& \langle z(t)-\varphi_n(u^*_n)(t),u^*_n(t)-u(t)\rangle +\langle
     f(t,u(t))-f(t,u^*_n(t)) ,u^*_n(t)-u(t)\rangle   \\
    &\le&\langle z(t)-g_n(u^*_n)(t),u^*_n(t)-u(t)\rangle +\langle
    g_n(u^*_n)(t)-\varphi_n(u^*_n)(t),u^*_n(t)-u(t)\rangle\\&+&\lambda_R(t)\|u^*_n(t)-u(t)\|^2.
 \end{eqnarray*}
 Integrating between $0$ and $t$, we obtain (since $u^*_n(0)=u(0)=u_0$)
 \begin{multline}\label{4.4}
   \frac 1 2\|u^*_n(t)-u(t)\|^2 \le \int_0^t\langle z(s)-g_n(u^*_n)(s),u^*_n(s)-u(s)\rangle ds\\
    +\int_0^t\langle g_n(u^*_n)(s)-\varphi_n(u^*_n)(s),u^*_n(s)-u(s)\rangle
    ds+\int_0^t\lambda_R(s)\|u^*_n(s)-u(s)\|^2ds.
 \end{multline} Since $\|g_n(u^*_n)-\varphi_n(u^*_n)\|_{\sigma} <\varepsilon_n$,
  then,  as in the proof of Theorem \ref{thm4.1} and by using Lemma \ref{lem2.3}, we have that $(g_n(u^*_n)-\varphi_n(u^*_n))$ converges
weakly in  $\mathrm{L}^1(I,H)$ to $0$,  and so, for all $t\in I$,
 \begin{equation}\label{4.5.1}
\lim_{n\to \infty} \int _0^t\langle
g_n(u^*_n)(s)-\varphi_n(u^*_n)(s),u^*(s)-u(s)\rangle ds =0.
\end{equation} On the other
hand, since $(g_n(u^*_n))$ and $(\varphi_n(u^*_n))$ are
 bounded in $\mathrm{L}^2(I,H)$ by $m(\cdot)$, and so in $\mathrm{L}^1(I,H)$, we have
 \begin{equation}\label{4.5}
\int _0^t\langle g_n(u^*_n)(s)-\varphi_n(u^*_n)(s),u^*_n(s)-u^*(s)\rangle
ds \le 2\|m\|_1\|u^*_n-u^*\|_C\underset{n\to \infty}\to 0.
\end{equation}
By \eqref{4.5.1} and \eqref{4.5} we get
$$ \lim_{n\to \infty} \int _0^t\langle g_n(u^*_n)(s)-\varphi_n(u^*_n)(s),u^*_n(s)-u(s)\rangle ds =0.$$
Also, we have by $(\mathcal{H}^6_F)$ and the fact that $z\in
S^2_{F(\cdot,u(\cdot))} $,
\begin{eqnarray*}
 && \int_0^t\langle z(s)-g_n(u^*_n)(s),u^*_n(s)-u(s)\rangle ds  \le \int_0^t\| z(s)-g_n(u^*_n)(s)\|\|u^*_n(s)-u(s)\|ds\\
   &\le& \int_0^t\Big(\varepsilon_n+d(z(s),F(s,u^*_n(s))\Big) \|u^*_n(s)-u(s)\|ds\\
   &\le&  \int_0^t\Big(\varepsilon_n+\mathcal{H}\big(F(s,u(s)),F(s,u^*_n(s)\big)\Big) \|u^*_n(s)-u(s)\|ds\\
   &\le&\int_0^t\Big(\varepsilon_n+k(s)\|u(s)-u^*_n(s)\|\Big) \|u^*_n(s)-u(s)\|ds.
\end{eqnarray*}
Consequently, letting $n\to \infty$ in \eqref{4.1}, we get
$$\frac 1 2\|u^*(t)-u(t)\|^2\le \int_0^t(k(s)+\lambda_R(s)) \|u^*(s)-u(s)\|^2ds.$$ Then, by Lemma \ref{lem2.4},
we obtain that  $u^* =u$. Consequently, $u^*_n\to
u$ in $\mathrm{C}(I, H)$ as $n\to \infty$ with
$(u^*_n)\subset S(\mathcal{P}_{ext(F)})$
 and this proves that $S(\mathcal{P}_{f,F})\subset \overline{S(\mathcal{P}_{f,ext(F)})}$.

 By using the same arguments in the proof of Theorem \ref{thm4.1}, we have that  $S(\mathcal{P}_{f,F})$ is closed
 in $\mathrm{C}(I, H)$.\\Since
 ${S(\mathcal{P}_{f,ext(F)})}\subset S(\mathcal{P}_{f,F})$, we conclude that $S(\mathcal{P}_{f,F})= \overline{S(\mathcal{P}_{f,ext(F)})}$,
 and this completes our proof. \end{proof}
\vskip2mm
 \textbf{Acknowledgements.} The authors thank professor  M.D.P. Monteiro Marques for some interesting
 discussions on the subject.

\end{document}